\documentclass{IEEEtran}

\usepackage{graphicx}
\usepackage{amsmath}
\usepackage{amsfonts}
\usepackage{amssymb}

\newtheorem{theorem}{Theorem}[section]

\newtheorem{definition}{Definition}[section]

\newtheorem{remark}{Remark}[section]

\newtheorem{assumption}{Assumption}[section]

\usepackage{atbegshi}
\AtBeginDocument{\AtBeginShipoutNext{\AtBeginShipoutDiscard}}
\begin{document}
\title{Stochastic Model Predictive Control, Iterated Function Systems, and Stability}
\thanks{Corresponding author V. Kungurtsev. V. K. and J. M. were supported by the OP VVV project
CZ.02.1.01/0.0/0.0/16\_019/0000765 ``Research Center for Informatics''. This work has received funding from the European Union’s Horizon Europe research and innovation program under grant agreement no. 101070568.}
\author{Vyacheslav Kungurtsev, Jakub Mare{\v c}ek, Robert Shorten}
\maketitle

\begin{abstract}                          
We present the observation that the process of stochastic model predictive control can be formulated in the framework of iterated function systems. The latter has a rich ergodic theory that can be applied to study the system's long-run behavior. We show how such a framework can be realized for specific problems and illustrate the required conditions for the application of relevant theoretical guarantees.
\end{abstract}

\section{Introduction}

Consider a generic iterative process with the following steps: a system takes a control input and behaves according to some stochastic dynamics; i.e., its output is noisy and governed by some (known or unknown) probability distribution. The control, in turn, is computed at each time step in order to solve a stochastic optimization problem that approximates these dynamics. This can be formalized as a map:
\begin{equation}\label{eq:mpcmap}
x(k)\to u(x(k),\eta_k) \to x(k+1,\xi_k)
\end{equation}
where the first $\to$ indicates a solution to a stochastic optimization problem given the input state $x(k)$,
\begin{equation}\label{eq:stochprog}
\begin{array}{l}    \min\limits_{u\in\mathcal{U}} \, \mathcal{R}\left[f(x(k+1,\xi),u,x(k))\right],\\ \text{such that }
    x(k+1,\xi) = g(u,x(k),\xi),\forall \xi \in \Xi
    \end{array}
\end{equation}
where $\xi$ is a random variable sampled from space $\Xi$ and $\mathcal{R}(\cdot)$ is a statistical aggregation operator, such as a risk measure, $f$ is a cost function, and $g$ describes some noisy dynamics.\newline 

In practice, typically~\eqref{eq:stochprog} cannot be solved exactly, but only by means of sample average approximation (SAA), in which Monte Carlo (or other stochastic discretization) samples of $\xi$ are taken and the optimization problem on the average is solved. SAA approximations satisfy the law of large numbers and are consistent (although biased) estimators (see, e.g.,~\cite{shapiro2014lectures}). Generically, this can be written as the first $\to$
in the schema~\eqref{eq:mpcmap} is a noisy operation, which we can consider to be subject to stochastic error $\eta$. We note that, alternatively, one can solve a more conservative variant of~\eqref{eq:stochprog}, such as a robust formulation (finding the optimal for the worst-case instance) or a distributionally-robust formulation (finding the optimal for the worst-case probability distribution among a set of possible ones). However, these methods would solve a different problem and are associated with their own advantages and drawbacks, which we do not consider in this note.\newline 

The second $\to$ in~\eqref{eq:mpcmap} corresponds to the stochastic realization of the subsequent state. Given the control calculated $u_k:=u(x(k),\eta_k)$ in iteration $k$, the next state satisfies the stochastic system equations $x(k+1,\xi_k) \sim g(u_k,x(k),\xi)$, where $\xi_k$ is the next realization in the stochastic process. Thus, with the distribution depending on the control $u_k$, the resulting output is another noisy function. This generic framework, although potentially modeling a variety of procedures, can be described as Stochastic Model Predictive Control (see, e.g.,~\cite{mesbah2016stochastic}), the iterative management of some physical process that is subject to random noise with known statistical properties.\newline 

\subsection{Motivation and Contribution}

There have been extensive algorithmic developments in solving the problem using approximations of the uncertainty using principled sample generation. However, in a significant thought-provoking article, Mayne~\cite{mayne2015robust} pointed out that stochastic MPC has been approached by the industry only hesitantly, due to important unresolved research questions. Specifically, the computational demands for even mildly nonlinear systems with uncertainty can often scale poorly with the degree of statistical confidence sought. Furthermore, studies of closed-loop stability have been limited and typically require a terminal constraint to be satisfied for all realizations, which translates to very conservative solutions of the corresponding optimization problems. The author challenges the scientific community to come up with schemes that perform extensive informative computations off-line, before the running of the control, and notions of stability that are more inclusive and coherent with uncertain dynamics. 
With this note, we wish to indicate that the framework of IFS can present a promising approach towards such a research program. \newline

Here, we propose a framework for reasoning about ergodic properties of stochastic model predictive control (SMPC). In particular, we place this process in the framework of an \emph{Iterated Function System} (IFS, cf. e.g.,~\cite{barnsley1985iterated,Diaconis1999}) to study the ergodic behavior of MPC problems of the form~\eqref{eq:mpcmap}. IFS describe a sequential probabilistic selection of maps to define a sequence of states of a system. IFS have been studied in various settings, with their long-run behavior studied through the lens of ergodic theory. This contrasts with standard notions of convergence of optimization problems and closed-loop stability, as considered in the traditional Model Predictive Control literature. These notions have proved challenging to extend to noisy systems due to their restrictiveness, and only recent results exist for stochastic problems, extending deterministic notions to the expectation~\cite{mcallister2022nonlinear}.  Ergodic theory has a rich set of conceptual and algorithmic tools as evidenced by the powerful monograph~\cite{meyn2012markov}. Thus, in this paper, we consider modeling this process as an IFS, and applying the relevant results and guarantees to the Stochastic MPC. In particular, we indicate what control problem's structural properties enable the application of theoretical results concerning the ergodicity. To the best of the authors' knowledge, the link has not yet been explored. Note that we do not introduce any new algorithms or solution procedures, but present the scaffolding for new means of analysis, which could provide understanding of the performance of existing algorithms. Of course, this in turn could provide insight into potentially effective techniques for novel procedures.\newline 

To maintain a generic formalism, we consider that $u$ and $x$ both live in some Polish space $\mathcal{X}$, and any use of a norm indicates its native norm. All functions, for example, $f$ and $g$ in~\eqref{eq:stochprog}, will be considered to live in $C^1(\mathcal{X})$, unless noted otherwise. $|\mathcal{A}|$ is the cardinality of a set $\mathcal{A}$.

\section{Iterated Function Systems}
\subsection{Background}

We begin by revisiting the notion of a \emph{state-dependent Iterated Function System} (IFS). This is a process wherein there exist a set of maps $\{F_i(x)\}$ and associated probabilities $p_i(x)$ where, at each step in the sequence, given the current state $x$, some index $i$ is chosen according to the probabilities $\{p_i(x)$\} and subsequently the map $F_i(x)$ is applied to generate the next iterate. A first foray into studying the properties of these maps is given in~\cite{barnsley1988invariant}. Although the literature has subsequently evolved considerably since publication of this paper, the paper remains a source of rich results for the community.\newline

\subsection{Discrete Controls and Exact Stochastic Programming Solutions}
Consider now the situation in which controls are discrete, i.e., there is a finite set $\mathcal{U}$ of possible inputs from which the control $u(k)$ must be chosen at each iteration $k$. The stochastic optimal control problem (OCP)  then amounts to taking the current $x(k)\in\mathcal{X}$, then computing an optimal $u(k)\in \mathcal{U}$ that minimizes the relevant probabilistic quantity, resulting in a noisy $x(k+1,\xi)\in\mathcal{X}$, or alternatively, computing an optimal mixed strategy $\{p_i\}$ of probabilities to implement $u_i\in\mathcal{U}$ with probability $p_i$.
(Note that since the expectation is a linear operator and we shall see that we require $f$ to be convex, we would expect the optimal control to be deterministic, i.e., $p_i=1$ for some $i$, if only the expected outcome is to be minimized or maximized. On the other hand, any higher moments or risk measures could make the mixed control optimal).\newline 

Formally, our goal is to solve for $p$ in the unit simplex $\bar\Delta$, such that the control $u_i$ is chosen among a finite set $\mathcal{U}$ with probability $p_i$:
\begin{equation}\label{eq:stochprogdisc}
\begin{array}{rl}
    \min\limits_{\{p_i\}\in \bar\Delta} & \mathcal{R}\left[f(x(k+1,\xi,i),u^p,x(k))\right],\\
    & u^p\sim \{u_i \text{ w.p. } p_i\},\\
    
    \text{s.t. } &
    x(k+1,\xi,i) = g(u_i,x(k),\xi),\forall \xi\in\Xi,\,i\in \text{supp} \{p_i\}
    \end{array}
\end{equation}
Note that we can then write $p_i(x(k))$ as depending on the previous state $x(k)$. The resulting state $x(k+1,\xi_k,i)$ depends on the control chosen $u_i\in \mathcal{U}$ according to $i\sim \{p_i\}$ and the physical realization of the noise $\xi_k$. Thus, by defining $S_i(\cdot)$ as the mapping from $x(k)$ to $x(k+1,\xi_k,i)$, we have shown that this procedure fits into the generic framework of a state-dependent IFS.\newline 

We shall now consider how one can apply the results on IFS, esp. \cite{tyran1997generic}, to conclude the long-run statistical properties of the behavior of~\eqref{eq:stochprogdisc}. To begin with, we must show that $p_i(x)$ as defined above satisfy the \emph{Dini condition}, which states that there exists a $\omega:[0,\infty)\to [0,\infty)$ continuous, non-decreasing and concave, with $\omega(0)=0$, such that $\int_0^1 \frac{\omega(t)}{t} dt<\infty$ and
$\sum\limits_{i=1}^{|\mathcal{U}|} \left| p_i(x)-p_i(y)\right|\le \omega(\rho(x,y))$ where $\rho$ is the metric on the underlying space. Indeed, by rewriting the problem as unconstrained,
\[
\begin{array}{l}
\min\limits_{\{p_i\}} \, R(\{p_i\}),\\
R(\{p_i\}) :=  \mathcal{R}\left[\sum_i p_i f(g(u_i,x(k),\xi),u_{i},x(k))\right]+1_{\bar\Delta}(\{p_i\})
\end{array}
\]
where $1_C(x)$ is the indicator of $x$ belonging to the set $C$, i.e., $1_C(x) = 0$ if $x\in C$ and $1_C(x) = \infty$ otherwise. Thus, in this case, $1_{\bar\Delta}(\{p_i\})$ enforces that $\{p_i\}$ lies in the unit simplex. Now, if we assume,\newline
\begin{assumption}
$\mathcal{R}(\{p_i\})$ is strongly convex with respect to $\{p_i\}$.\newline
\end{assumption}

we can now use~\cite[Theorem 4.1]{attouch1993quantitative}, considering $x(k)$ as the parameter of the problem, with the domain of $x(k)$ being any large enough compact set, to guarantee such a function $\omega(t)$ satisfying the Dini condition exists. Note, however, that if $\mathcal{R}=\mathbb{E}$ above, then the solution of the problem is clearly $p_i=1$ for $i$ such that $\mathbb{E}[f(g(p_i,x(k),\xi),u_{i},x(k))]$ is minimal, and this is a linear program, thus not strongly convex. This can be corrected simply by adding a regularization $\alpha \sum_i p_i^2$ to the objective. Second, we point out that this is a sufficient, but by no means necessary assumption.
It is our intention to open the field of analyzing SMPC with IFS, and results must by necessity begin with the most straightforward cases.

Finally, we use~\cite[Theorem 1]{tyran1997generic} to prove the ergodicity of the resulting IFS, which we state below,\newline

\begin{theorem}\cite[Theorem 1]{tyran1997generic}
\label{thm21}
Let $(S,p)$ be an iterated function system, i.e., there exist $S_i:\mathcal{X}\to \mathcal{X}$ for $i=1,...,N$ such that given $x$, with probability $p_i(x)$, the next state is defined by $S_i(x)$. If,
\begin{enumerate}
    \item There is a Dini function of $(S,p)$
    \item For every $i\in\{1,...,N\}$, $\inf\limits_{x\in \mathcal{X}} p_i(x) >0$  \label{cond2}
    \item The transformations $S_i:\mathcal{X}\to \mathcal{X}$ are $L(S_i)$-Lipschitzian for $i=1,...,N$ and there exists $\lambda_S$ such that,
    \[
    \sum\limits_{i=1}^N p_i(x) L(S_i) \le \lambda_S < 1\text{ for }x\in X,
    \] \label{cond3}
\end{enumerate}
then the system $(S,p)$ is ergodic, i.e., for the kernel operator $\mathcal{S}$ of the system, there exists a stationary distribution $\mu^*(x)$ such that for any initial $\mu_0(x)$, it holds 
\[
\lim\limits_{n\to\infty} \|\mathcal{S}^n \mu_0-\mu^*\|_{\mathcal{F}}
\]
where $\mathcal{F}$ indicates the Fortet-Mourier norm,
\[
\|\mu\|_{\mathcal{F}} = \sup\left\{\int f d\mu,\,f\in C(\mathcal{X}),\,\sup\limits_{x\in\mathcal{X}} |f(x)|\le 1 \right\}.
\]\newline
\end{theorem}

To apply Theorem \ref{thm21}, we must check the other two conditions. Condition \ref{cond2} implies that there is some $p_0$ such that for all possible states $x$, we have that $p_i(x)>p_0$. One sufficient condition for this to hold is that for all $i,j\in [N]$, we have some bound on the cost difference
$f(g(u_i,x,\xi),u_i,x)-f(g(u_j,x,\xi),u_j,x)$ that holds across $x\in\mathcal{X}$, possible control selections $i$ and noise $\xi\in \Xi$.\newline

Condition \ref{cond3} of Theorem \ref{thm21} requires Lipschitzianity. In particular, it must hold that for all maps $F_i(\cdot)$ from $x(k)$ to $x(k+1,\xi,i)$ are Lipschitzian with respect to $x(k)$, i.e., $g$ is Lipschitzian with constant $L_i$ with respect to the second argument. In addition, it must hold that,
\begin{equation}\label{eq:lipcontract}
\sum\limits_{i=1}^{|\mathcal{U}|} p_i(x) L_i < 1
\end{equation}
for all possible $x$, formally,
\[
\left\|g(u_i,x,\xi)-g(u_i,y,\xi)\right\|\le L_i\|x-y\|,\,a.e.\,\xi
\]
With this, we can now claim that the conditions of~\cite[Theorem 1]{tyran1997generic} hold. Thus, we can conclude that the IFS system defined by the repeated stochastic OCP is asymptotically stable. 

\section{Stochastic MPC Modeled as a IFS}
\subsection{Set Up}
Given the entire nested noise admixture of~\eqref{eq:mpcmap},
even the state-dependent IFS form~\cite{tyran1997generic} as considered above is not sufficiently expressive to adequately model the more general stochastic MPC process. 
For the more general case we consider solving~\eqref{eq:stochprogdisc} by first taking $J$ samples $\{\xi_j\}\sim \Xi$, and denoting this finite set as $\bar{\Xi}$. Then we minimize a sample average of the optimization objective, i.e., a SAA approximation,
\begin{equation}\label{eq:stochprogdiscsaa}
\begin{array}{rl}
    \min\limits_{u\in\mathcal{X}} & \mathcal{R}\left[f(x(k+1,\xi^p),u,x(k))\right],\\
    & \xi^p\sim \{\xi_j \text{ w.p. } 1/J\} \\
    \text{such that } &
    x(k+1,\xi_j) = g(u,x(k),\xi_j),\forall \xi_j\in \bar{\Xi}
    \end{array}
\end{equation}
Recall now the two sources of noise, when considered as a map from $x(k)$ to $x(k+1)$. First, $\xi_j$ themselves are sampled from $\xi_j\sim \Xi$. The sampling affects the outcome of solving the optimization problem, i.e., $u$ depends on the $J$ samples, $u(\{\xi_j\})$. Next, the state at time $k+1$ depends on the physical realization of $\Xi$, i.e., $x(k+1,\xi)$.
To model this, we must incorporate the notion of an  IFS~\cite[e.g.]{horbacz2001continuous,bielaczyc2007generic}, formally a pair $(S,p)$ with probability map $p(t,x)$ on state $x$ with parameter $t$ satisfying,
\[
\int_0^K p(t,x)dt = 1
\]
and Markov operator
\[
P_{(S,p)}\mu(A) = \int\limits_{\mathcal{X}} \int_0^K \mathbf{1}_A(S(t,x))p(t,x)dt\mu(dx)
\]
for $A\in\mathcal{B}(X)$, the Borel set on $\mathcal{X}$. Procedurally, given a state $x_k$, the probability density $p(\cdot,x_k)$ governs the realization of the continuously indexed map $S(t_k,x_k)$, which is itself deterministic.\newline 

To utilize the theoretical results associated with IFS, the process~\eqref{eq:mpcmap} must be appropriately linked to the underlying abstractions. Specifically, $p(t,x)$ must incorporate both the SAA noise $\eta$ and the output system noise $\xi$ into the parameter $t$. Then the map $S(t_k,x_k)$ corresponds to the output realization $x(k+1,\xi_k)$ for $t_k=(\eta_k,\xi_k)$. 

Now we introduce several notions from~\cite{horbacz2001continuous} associated with an IFS $(S,p)$ and its Markov kernel $P$.  
Recall that the \emph{dual} $U$ of $P$ is given by,
\[
\langle U f,\mu\rangle = \langle f,P\mu\rangle \text{ for }
f\in B(\mathcal{X}),\, \mu\in \mathcal{M}_{f}
\]
where $\mathcal{M}_{f}$ is the set of finite measures and $B(\mathcal{X})$ Borel measureable functions. 
The operator $P$ is called \emph{Feller} if $Uf\in C(\mathcal{X})$ for 
$f\in C(\mathcal{X})$ and \emph{nonexpansive} if $\|P\mu_1-P \mu_2\|_L
\le \|\mu_1-\mu_2\|_L$ for $\mu_1,\mu_2\in \mathcal{M}_f$. The dynamic properties of interest associated with this operator are notions of stability, convergence and ergodicity -- broadly speaking the limiting behavior of the probability distributions of the state. A measure $\mu\in \mathcal{M}_f$ is \emph{stationary} or \emph{invariant} if $P\mu=\mu$. The operator $P$ is \emph{asymptotically stable} if there exists a stationary distribution $\mu_*$ and constant $q>0$ such that, for any $\mu\in \mathcal{M}_f$
\[
\lim\limits_{n\to\infty} \|P^n \mu - \mu_*\|_q = 0
\]
Denote the limit points of the sequence of measures defining the process by,
\[
\omega(\mu) = \{\nu \in \mathcal{M}_f:\exists_{m_n,n\ge 1} m_n\to \infty\text{ and }P^{m_n}\mu\to \nu\}
\]
Let $\mathcal{C}_{\epsilon}$ be the family of all sets $C\in\mathcal{B}(\mathcal{X})$ for which there exists some finite cover of points with radius $\epsilon$, i.e., $\exists n$ and $\exists\{x_1,...,x_n\}\subset \mathcal{X}$ such that $C\subseteq \cup_{i=1}^n B(x_i,\epsilon)$.\newline

\begin{definition}
The operator $P$ is \emph{semi-concentrating} if for every $\epsilon >0$ there exists $C\in \mathcal{C}_{\epsilon}$ and $\alpha>0$ such that,
\begin{equation}\label{eq:semiconc}
    \lim\inf\limits_{n\to\infty} P^n \mu(C)> \alpha \text{ for }\mu\in\mathcal{M}_f\newline
\end{equation}
\end{definition}

Now let us consider explicitly the Markov operator
\begin{equation}\label{eq:moperator}
P_{(S,p)}\mu(A) = \int\limits_X \int_0^K \mathbf{1}_A(S(t,x))p(t,x)dt\mu(dx)
\end{equation}
for $A\in\mathcal{B}(X)$. Now if,
\begin{equation}\label{eq:lipcond1}
d(S(x,t),S(y,t))\le \lambda(x,t)d(x,y)
\end{equation}
with
\begin{equation}\label{eq:lipcond2}
\int_0^T \lambda(x,t) p(x,t)dt\le \gamma<1
\end{equation}
and,
\begin{equation}\label{eq:lipcond3}
\int_0^T|p(x,t)-p(y,t)|dt\le \theta d(x,y)
\end{equation}
with $\theta>0$ we have a stability result of the following form.\newline

\begin{theorem}\cite[Theorem 4.3]{horbacz2001continuous}
If $(S,p)$ satisfy conditions~\eqref{eq:lipcond1}-\eqref{eq:lipcond3}, then $P_{(S,p)}$ is semi-concentrating.\newline
\end{theorem}

Finally, an additional technical stopping-time condition provides a sufficient mechanism to ensure asymptotic stability for $P_{(S,p)}$.\newline

\begin{theorem}\cite[Theorem 4.3]{horbacz2001continuous}
Let $(S,p)$ satisfy the conditions~\eqref{eq:lipcond1}--\eqref{eq:lipcond3}. In addition, assume that there exists an $\gamma>0$ such that for all $x\in \mathcal{X}$, there exists a time $\tau_x\in[0,T]$ satisfying,
\begin{equation}\label{eq:taux}
p(x,t)=0,\,\text{for }0\le t<\tau_x\text{ and }p(x,t)>\gamma\text{ for }\tau_x\le t\le T
\end{equation}
and $p(x,\cdot):[\tau_x,T]\to \mathbb{R}_+$ is continuous. If, in addition, $\sup_{x\in X} \tau_x<T$ then $P_{(S,p)}$ is asymptotically stable.\newline
\end{theorem}

\subsection{Key Observations and Open Questions}

Consider now two possible initial states at $k$, $x$ and $y$ and a potential set of realizations $t=(\eta,\xi)$. Solving the optimization problem, subject to SAA noise, defines $p(t,x(k))$ as the distribution of chosen $u(x(k),\eta_k)$ with the noise defining the homotopy with respect to $t$. This in turn induces the map $S(t,x)$ as defined by $g(u(x(k),\eta),x(k),\xi_k)$ once $t_k=(\eta_k,\xi_k)$ has been chosen.\newline 

To consider the assumptions in the corresponding state dependent IFS theory, if for all SAA realizations $\eta$, the optimization problem~\eqref{eq:stochprogdiscsaa} is Lipschitz stable with respect to the input $x(k)$ which holds, e.g., if the map $\mathcal{R}(x(k+1,\xi^p(\eta),u),u,x(k))$ is strongly convex with respect to $u$, where we now write the subsequent state as a function of $u$ (the reduced problem), then clearly this will hold globally. Otherwise, for non-convex objectives with local minimizers that satisfy second-order sufficient conditions for optimality, this condition holds locally.\newline 

We recall second-order sufficient conditions for optimality (e.g.,~\cite{nocedal2006numerical}).
\begin{definition}
The \emph{second order sufficient conditions for optimality} conditions hold at $u$ if for all $\Delta u\neq 0$ it holds that,
\[
\langle\Delta u,\nabla^2_{uu}\mathcal{R}(x(k+1,\xi^p(\eta),u),u,x(k)) \Delta u\rangle > 0
\]
\end{definition}
See the results of~\cite{robinson1982generalized} for upper Lipschitz continuity of the optimal solution $u$ as a function its parameters, which in this case corresponds to $x(k)$.\newline

\begin{remark}
In many cases, $u$ is required to exist in some compact bounded set $\mathcal{U}$. This introduces the necessity to consider active sets in the formulation of the second-order optimality conditions, which add additional notation without additional insight here. Note, however, that in~\cite{bielaczyc2007generic} it is shown that if we are constrained to a compact convex set, the invariant measure has Hausdorff dimension zero.\newline 
\end{remark}

Subsequently, if also $g(u,x,\xi)$ being Lipschitz stable as a function of $x$ as well for all $\xi$ and $u$ then we have achieved sufficient conditions for there being some $L$ such that $|p(x,t)-p(y,t)|\le Ld(x,y)$.\newline 

\textbf{Open Problem:} In order to utilize this approach, the techniques of upper Lipschitz continuity subject to perturbations (e.g., in the comprehensive monograph~\cite{bonnans2013perturbation}) need to become \emph{quantitative}, to obtain estimates of the moduli of continuity or appropriate scaling metrics.\newline 

Let us now turn to the other condition, given by~\eqref{eq:taux}. In the context of our stochastic OCP, this implies that certain distributions of noise $t=(\eta,\xi)$ are inaccessible for some states $x(k)$. Since $\eta$ can be regarded as exogenous, or state independent, as determined by the samples generating the SAA, it must imply that for certain $x(k)$, there are $u$ that do not lie in the support of the distribution of solutions of the SAA problem across realizations $\eta$. Such a question cannot be answered without distributional information with respect to the noise structure of the dynamic process.\newline

\textbf{Open Problem(s):} Taking into context the specific problem-dependent distributional information $x(k+1)\sim g(u,x(k))$, characterize the conditions of finite support of $u^*$ as a solution of the optimization problem defined by SAA sampling.\newline 

\section{A Numerical Illustration}
We performed a synthetic simulation to illustrate the ergodicity of states in stochastic model predictive control. We consider the four-state system,
\[
x_{k+1} = (A+\Xi) x_{k}+B u
\]
where $\Xi$ is additive noise. The objective function is the standard MPC tracking objective with regularization:
\[
f(x,u)=\mathbb{E}\left[(x-z)^T Q (x-z)\right]+u^T R u\newline
\]

We generated the problem as follows: 
\begin{enumerate}
    \item To encourage contractive dynamics, we took $\Lambda_A=\text{diag}\begin{pmatrix} 1/5 & 1/8 & 1/10 & 1/12\end{pmatrix}$ generated a random orthonormal eigenbasis $V$ and let $A=V^T \Lambda_A V$.
    \item The tracking and regularization matrices $Q$ and $R$ were similarly made to be positive definite, with $\Lambda_Q=\text{diag}\begin{pmatrix} 5 & 6 & 9 & 15\end{pmatrix}$ and $\Lambda_R=\text{diag}\begin{pmatrix} 0.5 & 2 & 1 & 1.5\end{pmatrix}$
    \item The control matrix $B$ has entires drawn from a uniform random variable over $[0,1]$, ensuring w.p. one full-rank.
    \item The target matrix $z$ is a four component vector drawn randomly uniformly from $[0,1]$.
    \item The perturbation matrix is of the form:
    \[
    \Xi = \begin{pmatrix} 0 & \xi_1 & 0 & 0 \\ 0 & 0 & 0 & 0 \\ 0 & 0 & \xi_2 & 0 \\ 0 & 0 & 0 & 0 \end{pmatrix}
    \]
    with $\xi_1,\xi_2\sim [-0.005,0.005]$
    \item We ran 20 trials (i.e., twenty such generations described) of 10000 iterations of stochastic MPC with $N=100$ samples for each SAA problem. Each output $x_{k+1}$ was computed from the system matrix with a random draw in $\Xi$.\newline 
\end{enumerate}

It can be seen that, since the process is linear, with $\mathbb{E}[A+\Xi]=A$, we can compute the exact control as the solution to,
\[
(R+B^T Q B)u = -B^T Q (A x_k-z)
\]
Thus, the state map satisfies,
\[
S(x,\Xi,\eta) = (A+\Xi) x-B\left((R+B^T Q B)^{-1}Q (Ax-z)+\eta\right) 
\]
where $\eta$ is the perturbation associated with solving the inexact SAA problem. For initial states $x$ and $y$ and given $\eta$ we have,
\begin{align*}
\|S(x,\Xi,\eta)-S(y,\Xi,\eta)\| \quad \quad \quad \quad \quad \quad \quad \quad \quad \quad \quad \quad \quad \\ \le 
\|A+\Xi\|\|x-y\|+\|(R+B^T Q B)^{-1}Q A\|\|x-y\|,
\end{align*}
implying that a sufficient condition for the contraction is:
\[
\|A+\Xi\|+\|(R+B^T Q B)^{-1}Q A\|<1\,\forall \Xi\sim \Xi(\xi_1,\xi_2) 
\]
This is equivalent to a stable dynamic system for every possible stochastic realization. Thus the fact that it agrees with ergdocitiy of the long-run stochastic behavior is indicative of a potential sea of relationships between the aggregation of point-realization dynamics and statistically agglomorated dynamics of the system.\newline

\begin{figure}[!h]
\centering
\includegraphics[scale=0.47]{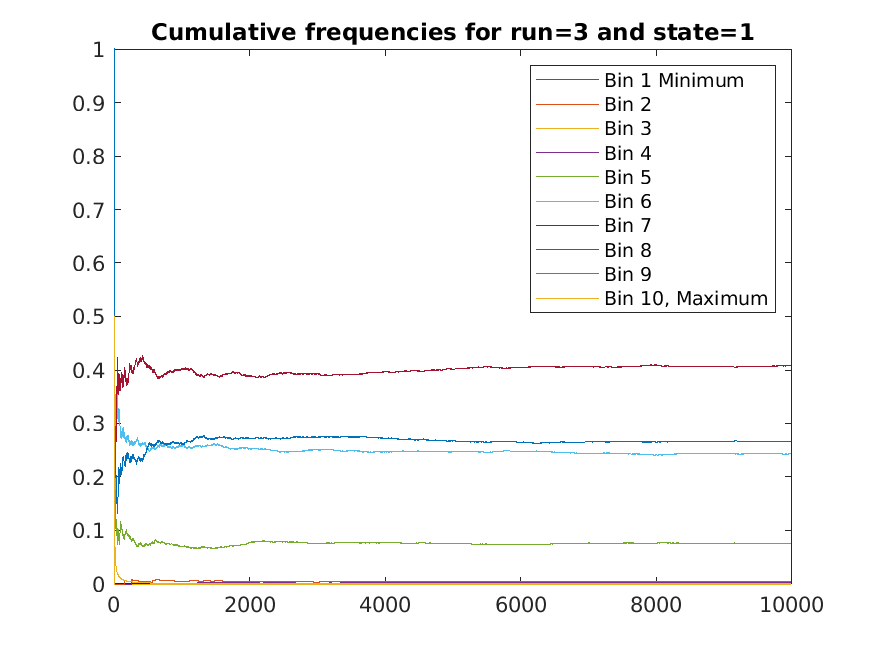}
\includegraphics[scale=0.47]{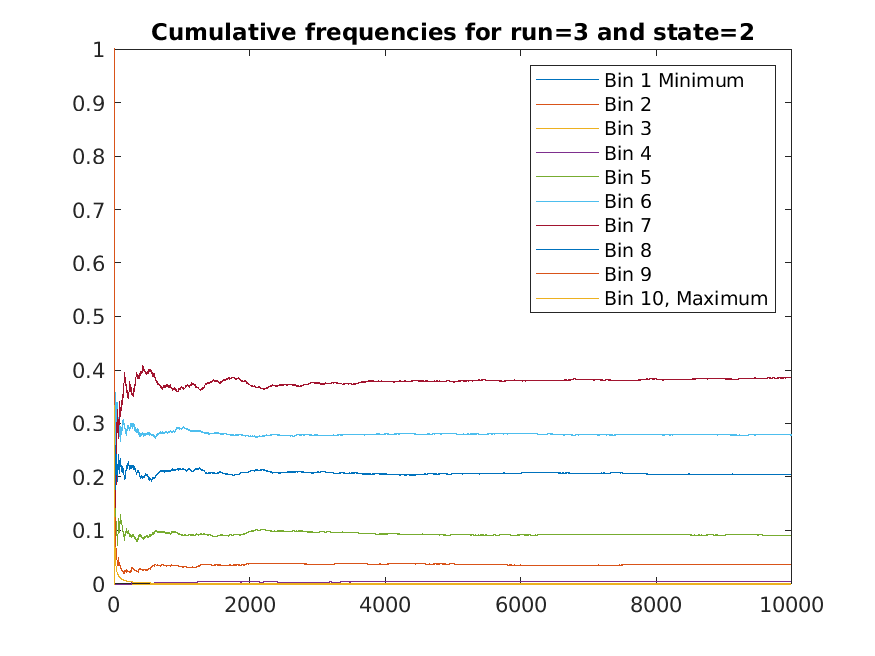}
\includegraphics[scale=0.47]{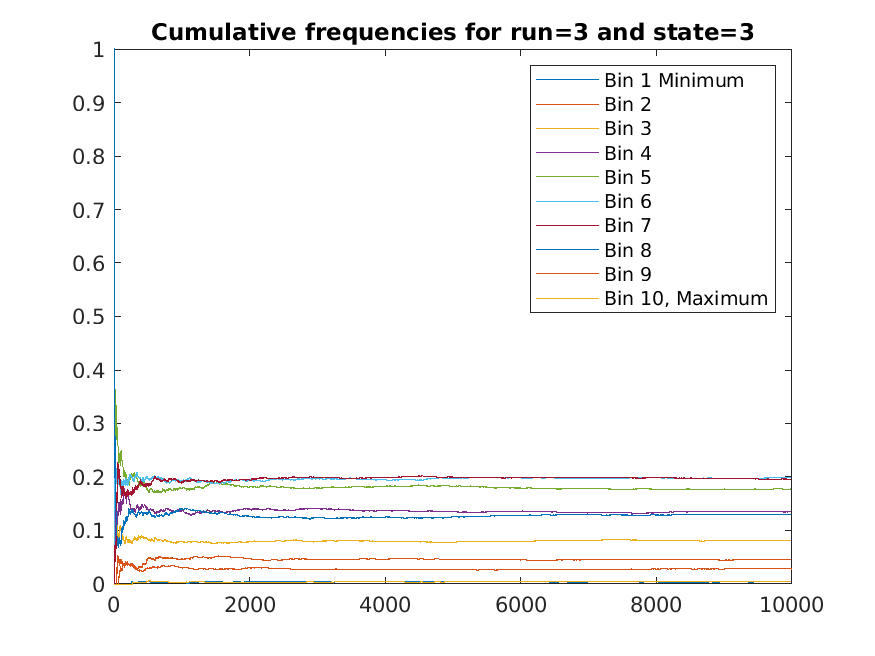}
\includegraphics[scale=0.47]{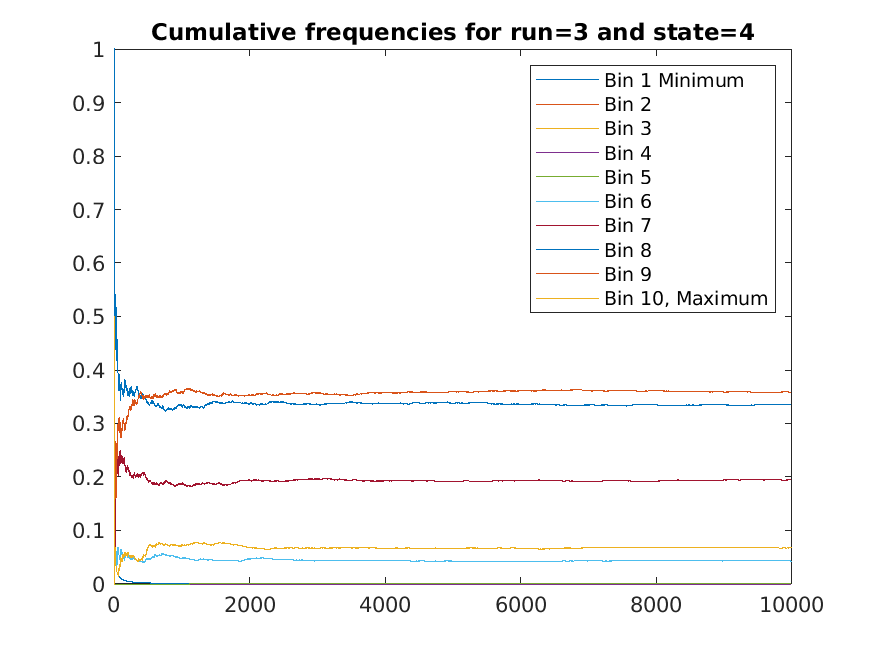}
\caption{Distributions across the bins for each state in the synthetic stochastic MPC run.} \label{fig:distributions}
\end{figure}

We computed the cumulative empirical distributions of the states as follows: for all four states, we took the minimum and maximum values over ten thousand iterations and defined a set of equally spaced discretizations as a histogram-bin, counting the proportion of times each state appears in each bin. Figure~\ref{fig:distributions} shows a representative run. It can be seen that the empirical distributions stabilize over the long run, suggesting the ergodicity of the process.\newline 

Furthermore, it can be seen that more restrictive notions of stability are uninformative in this case. In particular, one can see that although the frequencies across the bins stabilize, the system still traverses a large state space. This can be seen because the bins include the range of the transient dynamics and the asymptotic frequencies are non-zero for up to seven of the bins. Thus, any notion that relies on localization, i.e., that the trajectory asymptotically approaches remains in some comparatively bounded region, would either fail or be vacuous, as both the asymptotic and initial transient dynamics traverse a comparatively similarly large portion of the state space. On the other hand, the relative frequencies with respect to how much time it spends across components of the region does clearly stabilize. Thus, probabilistic notions of asymptotic occupation and convergence, rather than localized notions, are already seen to be more appropriate for this simple toy example.

\section{Conclusion}
The framework of IFS presents a rich and powerful set of tools in the analysis of limiting statistics of iterative processes. Stochastic MPC can be formulated in this framework; under appropriate conditions, important results can be proven regarding its behavior. These results depend on the conditions that can be studied when one has accurate (either a priori, or data driven) distributional information on the process. This, in turn, suggests a comprehensive program of applied analysis of such problems.

\bibliographystyle{plain}
\bibliography{refs}
\end{document}